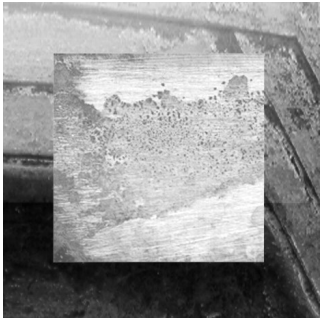

**Protopapas Mattheos**
University of Rome "La Sapienza"
**Kosmatopoulos Elias**
Technical University of Crete

# Simulation and Use of Heuristics for Peripheral Economic Policy

**Abstract.** Recent trends in Agent Computational Economics research, envelop a government agent in the model of the economy, whose decisions are based on learning algorithms. In this paper we try to evaluate the performance of simulated annealing in this context, by considering a model proposed earlier in the literature, which has modeled an artificial economy consisting of geographically dispersed companies modeled as agents, that try to maximize their profit, which is yielded by selling an homogenous product in different cities, with different travel costs. The authors have used an evolutionary algorithm there, for modeling the agents' decision process. Our extension introduces a government agent that tries to affect supply and demand by different taxation coefficients in the different markets, in order to equate the quantities sold in each city. We have studied the situation that occurs when a simulated annealing algorithm and a simple search algorithm is used as the government's learning algorithm, and we have evaluated the comparative performance of the two.
**Keywords:** Agent Computational Economics, Simulated Annealing, Peripheral Economic Policy.

## 1. INTRODUCTION

Agent Computational Economics (ACE) concerns the computational study of economies modeled as evolving systems of autonomous interacting agents (Testfatsion 2001). Various ACE models, concern the modeling of simple, artificial economies, in order to study the complex inter-correlations between economic entities (agents) and the macroeconomic phenomena they give rise to.

Ichibuchi et al. (2001), have proposed one such model. 100 geographically dispersed firms, modeled as agents, existed on a 100x100 grid, and were trying to maximize their profit by selling one unit of a homogeneous product, in one of the 5 markets that also existed on the grid. The agents used an evolutionary algorithm, in order to learn how to choose the correct market they would sell to. Ichibuchi et al. (2001) displayed, by means of simulating the model, that the agents would indeed maximize their profit.

We have extended the model, by introducing taxes, imposed by a hypothetical government. Taxation changes the equilibrium in the economy, and therefore allows the government to affect it, in way that serves its macroeconomic goals. A planner can simulate the above model for different sets of tax coefficients and use the observed outcome in order to make more thoughtful decisions about taxation matters.

We have pursued this goal even further, by introducing an optimization scheme based at the algorithm of simulated annealing. The algorithm simulates the model of the economy, for an adequate number of times, in order to achieve good estimates for the statistic that gives the objective function's value. Then it changes the tax coefficients, in the way proposed by the simulated annealing algorithm. Note that the choice of the objective function is subjective to the personal preferences of the planner. Here we have pursued the simple goal of trying to equate the quantities sold in each market.

In order to evaluate the efficiency of the simulated annealing algorithm, we have also used a simple stochastic search for the same problem. The difference in the values of the objective function, acquired by the two algorithms is evident.

## 2. THE MODEL

We, just like Ichibuchi et al. (2001) assume there are 100 firms and 5 markets in a geographical grid of 100x100 tiles. Firms produce one unit of a homogeneous product each round, and sell it in one of the 5 markets, in order to maximize their profit. Since agents produce one unit of product we can write

$$\sum_{j=1}^{5} x_{ij}^t = 1 \qquad (1)$$

$x_{ij}^t$ denotes the quantity of product agent *i* sells in market *j*, at time step *t*. The quantity the agent sells in the market he picks for a given round is 1, and the quantity he sells in all the other markets is 0.

We assume that the demand function is linear i.e. that the price in a market is a linear function of the sum of quantities sold in that market.

$$p^t_j = a_j - b_j X^t_j \tag{2}$$

$a_j$ and $b_j$ are exogenous coefficients, and $X^t_j$ is the sum of quantities sold in market $j$ at time step $t$ i.e.

$$X^t_j = \sum_{i=1}^{100} x^t_{ij} \tag{3}$$

Ichibuchi et al. (2001) have set $a_j = 200$ and $b_j = 3$. These are the values we are going to use as well.

Transportation costs are also assumed to be a linear function of the Euclidean distance between an agent and the market.

$$c_{ij} = c d_{ij} \tag{4}$$

$c_{ij}$ is the cost suffered by agent $i$ when he sells his product in market $j$, and $d_{ij}$ is the Euclidean distance between the two. We also follow Ichibuchi et al. on this one, so we set $c=1$.

The equation for agent's $i$ profit, follows directly from the inverse demand (2) and the cost (4) functions. Taking into account, that each firm sells exactly one unit of product each round, we have

$$r^t_{ij} = p^t_j - c d_{ij} \tag{5}$$

The evolutionary algorithm that determines the agents' action is, as introduced by Ichibuchi et al (2001), based on "mimicking" and "mutation". All agents make the choices simultaneously, and their profit determines the fitness value of their chosen strategies. If we denote $r(s^t_{ij})$ agent's $i$ profit, based on the strategy $s_{ij}$ he picks on round $t$, then this profit will be used for determining the agent's choice for round $t+1$.

"Mimicking" has the following meaning: In any given round, an agent "decides" if he is going to change his strategy, with probability $P_r$. If this is the case, he picks one of the strategies employed on the previous round, by one of his $N$ neighbors. These $N$ neighbors are the firms that are located closer to him on the grid, including itself. The probability of choosing a strategy is proportional to the profit that strategy yielded, for the firm that employed it, in the previous round. More precisely is given by

$$P_k = \frac{r_k - r_{\min}}{\sum_{k \in N(i)} r_k - r_{\min}} \tag{6}$$

$P_k$ is the probability the agent will choose the quantity his neighbor $k$ chose on the previous round, $r_k$ is the profit agent $k$ got from the use of this strategy, and $r_{min}$ is the minimum payoff amongst the neighbors. The value in the denominator ensures that probabilities add up to 1.

"Mutation" introduces a random change in an agent's strategy, happening with probability $P_m$, independent from $P_r$. Summing up, we can describe the evolutionary algorithm proposed by Ichibuchi et al. (2001) as follows:

1. Each agent's strategy for $t=1$ is chosen randomly.

2. Agents use their strategies at every time step. Then their profit is calculated, according to (5).

3. At every time step, any agent can change the strategy that he used at the previous step, with probability $P_r$. If this is the case, he imitates one of his $N$ neighbor's (the agents that are located the closest at him, on the grid, including itself). The probabilities for picking any of these strategies are given by (6).

4. Each agent can change his strategy, to a completely random one, with probability $p_m$. If he does, he picks one of all the feasible strategies with equal probability.

5. Loop back to step 2, until a termination criterion holds.

We have expanded the model, by introducing a government, that tries to affect the market state, by imposing taxes and giving beneficial payments to the firms. There are two kinds of those taxes (payments): taxes (payments) based on the price prevailing in a market

$$T_1 = t_j p_j^t \tag{7}$$

and taxes (payments) of a fixed amount, that can be used to "lure" firms in a given market, even more effectively.

$$T_2 = a_j \tag{8}$$

$t_j$ and $a_j$ can be either positive or negative. They are positive in the case of beneficial payments and negative in the case of taxes.

Hence, a firm's profit is no longer provided by (5), but is instead given by

$$r_{ij}^t = (1+t) p_j^t + a_j - c d_{ij} \tag{9}$$

since it's revenue is no longer equal to the price in the market (recall that the quantity it sells is always equal to one), but instead to the price, plus the beneficial payments minus the taxes. Replacing (1) by (9) causes –along with the minor change in the evolutionary algorithm- the equilibrium state, concerning both the price and total quantity in each market, to change. The government decisions, when it comes to establishing the values for the tax coefficients, have an impact on the economy.

Therefore it is important to try to select these coefficients in a way, that maximizes (or tends to maximize) a given macroeconomic criterion. This criterion is subject to choice of each macroeconomic planner. In this implementation, we assume the government is trying to minimize the deviations on the quantities sold in each market, given by

$$\frac{1}{900\sqrt{4}} \sum_{t=100}^{900} \sqrt{\sum_{j=1}^{5} \left(Q_j^t - \bar{Q}_j^t\right)^2} \quad (10)$$

The above is the average of the standard deviations of the quantities in the 5 markets for the time steps starting at $t=100$ and ending at the end of the simulation time $t=1000$. We have used the unbiased estimator for the variance

$$\hat{s}^2 = \frac{1}{n-1} \sum_{i=1}^{n} \left(x_i - \bar{x}\right)^2 \quad \text{(Koliva and Bora 1995)}$$

One could also use the standard estimator

$$s^2 = \frac{1}{n} \sum_{i=1}^{n} \left(x_i - \bar{x}\right)^2$$

or even neglect the constant coefficient of (10) entirely. All these formulations are equivalent, in a minimization problem, since the vector of tax coefficients that minimizes (10), will also minimize any function of the form

$$A \sum_{t=100}^{900} \sqrt{\sum_{j=1}^{5} \left(Q_j^t - \bar{Q}_j^t\right)^2} \quad (11)$$

Noting that we have 100 firms, each one producing one unit of good and 5 markets, we can also say that

$$\bar{q}_j^t = 20 \quad (12)$$

Since the evolutionary algorithm used for choosing the strategies for the agents, is stochastic in nature, a given set of tax coefficients would yield different values for (10), every time we run the simulation. This means that the application of meta-heuristic algorithms, in order to optimize the choice of tax coefficients, would be problematic, since these algorithms assume a deterministic objective function to optimize.

We could solve this problem by taking a good estimate of the expected value of the stochastic objective function. A good estimate is needed, in order to be as positive as possible, that the true value of the mean of the objective function is indeed better that the values given by other values of tax coefficients.

We are going to use simulated annealing meta-heuristic (Kirckpatrick et al. 1983) to find an approximation for the global optimum, of the minimization problem of (10), for a set of acceptable values for the tax coefficients $t_j, a_j$. The algorithm, given in pseudo-code, is as follows:

1. Pick initial solution $s_0$ and initial temperature $t_0 > 0$
2. Choose a sequence for updating the temperature
$$t_{+1} = a(t)$$
3. Repeat the following outer loop, until a termination criterion happens
4. Repeat the following inner loop, until a termination criterion happens
5. Choose a random solution from a neighborhood of the initial solution $s \in N(s_0)$
6. $d = f(s) - f(s_0)$
7. If $d < 0$ update initial solution $s \rightarrow s_0$
8. else pick a random $x \in (0,1)$
9. if $x < \exp\left(-\dfrac{d}{t}\right)$
10. then update solution $s \rightarrow s_0$
11. end of inner loop
12. set $t \leftarrow a(t)$
13. end of outer loop.

This approach can be regarded as a variant of the well-known heuristic technique of local search (a variant of which we will also employ here), in which a subset of the feasible solutions is explored by repeatedly moving from the current solution to a neighboring solution. For a minimization problem the traditional forms of local search employ a descent strategy, in which the search always moves in a direction of improvement. However such a strategy often results in convergence to a local rather than a global optimum.

On the other hand, simulating annealing allows some uphill moves in a controlled manner, and offers a way of overcoming this problem. The probability of moving from one state to another can be represented in matrix form, and for constant temperature, it only depends on the initial and final states.

This type of transition matrix gives rise to what is known as a *homogeneous Markov chain* (Van Laarhoven and Aarts 1988). As long as it's possible to find a sequence of exchanges which will transform any solution in any other with non-zero probability, the process converges towards a stationary distribution which is independent from the starting solution. As the temperature tends to zero, so the form of the stationary distribution tends to a uniform distribution over the set of optimal states or solutions (Van Laarhoven and Aarts 1988).

The parameters needed for the algorithm are the initial temperature and the sequence of temperatures in general, and the size (and structure) of the neighborhood of the initial solution on step 5. We also need to set the stopping criteria for the inner and outer loops.

In order to have a comparison measure for the algorithm, we will also use a simple (stochastic) local search algorithm in our problem:

1. Pick initial solution $s_0$
2. Loop until termination criterion occurs
3. Pick random solution in a neighborhood of the current
   $s \in N(s_0)$
4. If it is better than the current solution, i.e.
   $f(s) < f(s_0)$
   update the current solution
   $s \to s_0$
5. end of loop.

So the entire optimization scheme we are going to use, can be described as follows:
   1. The meta-heuristic algorithm sets the values of the tax coefficients
   2. The simulation is executed for $N_{sim}$ number of times.
   3. An estimator for the expected value of (10) is derived.
   4. The estimator for the objective function's value is used, in order to select the next solution (values of tax coefficients), according to the heuristic algorithm.
   5. Loop the above steps (1)-(4) for a given number of times.

As we mentioned above, a $N_{sim}$ number of independent repetitions of the simulation is needed, in order to acquire a "good" estimate, for the expected value of the objective function. The number of repetitions required to estimate this with a given degree of accuracy, depends on the variance of the estimator. In a sample of $n$ independent measurements, the real value of the mean will reside in the interval

$$\left(\bar{x} - z_{\frac{a}{2}} \frac{s}{\sqrt{n}}, \bar{x} + z_{\frac{a}{2}} \frac{s}{\sqrt{n}}\right) \text{ for } n > 30 \qquad (13)$$

with probability $1-a$ (Koliva and Bora 1995). We need therefore an estimate for the s.t.d. as well.

## 3. SIMULATION SETTINGS

For Ichibuchi et al (2001) evolutionary algorithm, we have used the same parameters, as the ones shown to be optimal in their study. So we have

$$N = 4 \quad P_r = 0.5 \quad P_m = 0 \qquad (14)$$

We have tested the algorithms in two different artificial economies. Placement of firms and markets on the grids of the two economies is random.

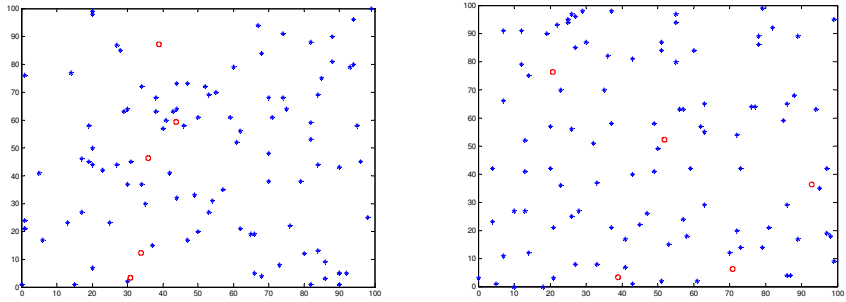

We've assumed the following feasible sets for the values of the tax coefficients:

$$t_j \in [-.25, .25] \quad a_j \in [-50, 50] \qquad (15)$$

In order to determine the standard deviation estimator in (13) we have taken several samples in each market, each sample based on different values for the tax coefficients. The estimated s.t.d. of the expected value of the objective function, never reached 4.5 Setting $s = 4.5$ in (13) implies that

$$N_{sim} = 10000 \qquad (16)$$

yields, with probability ¿.95 , that the observed value of the objected function, would be no different from the actual expected value by .05

$$p(|\bar{X} - \mu| > 0.05) < 0.05 \tag{17}$$

One could improve the absolute error even more, by setting a greater value for $N_{sim}$, provided he has access to the necessary computing power. By accepting (16), about 4 days of computation in a Pentium 4 are required, for a single run of the simulated annealing algorithm.

When it comes to the neighborhood structure, we have used the same structure for both algorithms. A random vector of new coefficients is picked, in every step, such that

$$t_j^{+1} \in \left[ t_j^0 - .02, t_j^0 + .02 \right] \quad , \qquad a_j^{+1} \in \left[ a_j^0 - 5, a_j^0 + 5 \right] \tag{18}$$

in order to construct the new candidate solution for the given step.

For the rest of simulated annealing parameters, we have used

$$t_{+1} = .8t \tag{19}$$

for the updating sequence of the temperature. Values in the range of .8 to .99 are usually used for this linear sequence. Since we lacked the computational power, we had used the minimum of these values. Using a greater value, would probably improve results even more. Values greater than 10 for the initial temperature, haven't improved the results, while a termination criterion of

$$t_f < .001 \tag{20}$$

is enough since, given the limited precision of any computer implementation, as temperature approaches zero, the small but positive probability of accepting any uphill move will be indistinguishable from zero. Finally we have used a fixed number of 10 internal iterations. The above parameters imply that the total number of steps in any execution of the algorithm will be 210.

We have used the same neighborhood structure for the stochastic local search algorithm, as we've already said. We have also used a termination rule of exactly 200 iterations, in order to make results comparable with these of the simulated annealing.

## 4. RESULTS

We have first taken a sample of 30 runs of the simulation (the evolutionary algorithm) for 3 different sets of tax coefficients, for each one of the markets.

Table 4.1. Sets of Coefficients $t_1, \ldots, t_5, j_1, \ldots, j_5$

| 0 | 0   | 0   | 0   | 0  | 0 | 0  | 0  | 0  | 0  |
|---|-----|-----|-----|----|---|----|----|----|----|
| 0 | .2  | .15 | .6  | .1 | 0 | 10 | 20 | 15 | 0  |
| 0 | .1  | .15 | .05 | .1 | 0 | 10 | 20 | 15 | 50 |

Table 4.2. Observed values of the objective function in each sample.

| Market 1 | | | Market 2 | | |
|---|---|---|---|---|---|
| 3.2755 | 6.0784 | 6.7937 | 4.5354 | 10.7729 | 2.7578 |
| 4.7592 | 5.9891 | 6.3925 | 3.4531 | 7.7866 | 3.4459 |
| 3.3443 | 9.4113 | 7.3426 | 5.6203 | 8.6590 | 6.7661 |
| 3.7415 | 5.3049 | 5.5322 | 3.4730 | 7.7729 | 3.9206 |
| 4.1054 | 7.6613 | 7.3578 | 3.9682 | 8.3326 | 7.8242 |
| 3.8024 | 5.9071 | 5.4336 | 2.8711 | 12.1485 | 3.5626 |
| 4.3517 | 5.9493 | 6.2533 | 4.6562 | 8.0667 | 4.1492 |
| 5.1756 | 2.8445 | 8.2043 | 2.2026 | 12.0682 | 5.8868 |
| 3.9144 | 8.3941 | 5.9939 | 4.9655 | 8.3363 | 3.1660 |
| 3.2440 | 8.6456 | 5.0737 | 4.0912 | 10.7316 | 3.9537 |
| 3.5894 | 7.2586 | 6.8270 | 3.5702 | 11.8253 | 3.3435 |
| 4.8005 | 6.1737 | 8.1118 | 8.7940 | 6.6583 | 3.7083 |
| 3.4899 | 6.1850 | 5.1340 | 4.9512 | 7.8316 | 2.5013 |
| 4.1390 | 6.2460 | 7.6728 | 4.7675 | 11.8774 | 5.3347 |
| 4.2502 | 5.6287 | 4.8842 | 4.0566 | 11.8673 | 5.4104 |
| 4.2622 | 7.6761 | 7.0582 | 3.9684 | 6.6500 | 7.5940 |
| 5.7649 | 5.8351 | 5.6151 | 5.0871 | 8.3826 | 3.9600 |
| 5.0027 | 6.3808 | 4.5798 | 3.8202 | 8.8585 | 4.6192 |
| 6.7380 | 5.1408 | 5.5507 | 6.0489 | 6.6547 | 4.0163 |
| 4.8580 | 5.3611 | 3.8015 | 6.0965 | 10.9842 | 5.0207 |
| 3.0874 | 7.5285 | 6.2968 | 3.7025 | 8.7010 | 4.8366 |
| 3.8201 | 5.7912 | 6.0111 | 8.6475 | 8.6815 | 4.1839 |
| 3.4217 | 4.7998 | 6.4525 | 4.5835 | 8.6451 | 4.3282 |
| 2.4953 | 6.9955 | 6.2964 | 5.6866 | 8.6730 | 3.1542 |
| 2.6535 | 5.5954 | 4.2684 | 6.7585 | 5.7041 | 3.9683 |
| 4.6259 | 8.7600 | 6.3606 | 4.0161 | 11.4358 | 4.1234 |
| 2.9554 | 5.6327 | 6.2271 | 4.9915 | 10.7887 | 3.6331 |
| 3.6533 | 5.7146 | 6.9510 | 2.9042 | 8.2262 | 2.7741 |
| 3.6514 | 6.1186 | 7.2337 | 3.8027 | 10.0296 | 4.0025 |
| 3.3503 | 4.8884 | 5.0169 | 5.3298 | 8.0567 | 3.6825 |

The qualitative behaviour of the evolutionary algorithm is the same, in the cases of no taxes (Ichibuchi et al) and our implementation, as shown at the two indicative executions, that follow.

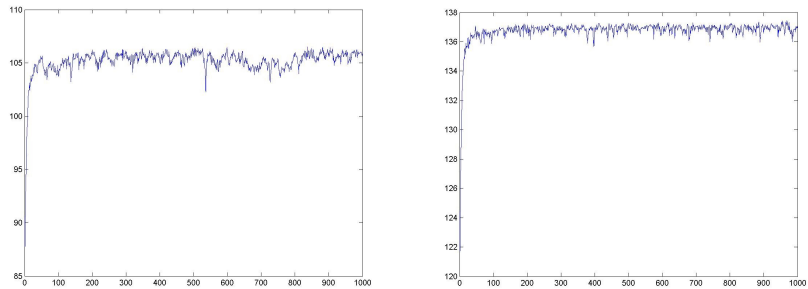

Fig 4.1. Mean profit for an execution with no taxes (left) and taxes different from 0 (right)

In order to be able to compare the performance of the two algorithms, we have to take samples of at least 30 executions of each. Since we do not know the exact distributions of the observed fitness (objective function) values, we need this number of executions, in order to approximate the unknown distributions by the normal distribution (Koliva and Bora 1995). The values that the two algorithms have converged to, are shown in tables that follow:

Table 4.3.1 Observed values of the objective function for 30 executions of simulated annealing, in the first economy

| 2.7163 | 2.5474 | 2.6277 | 2.2418 | 2.7115 | 2.1313 |
|---|---|---|---|---|---|
| 4.0827 | 2.8485 | 2.4052 | 2.4604 | 6.575  | 3.1888 |
| 2.4102 | 2.8882 | 2.4225 | 2.8714 | 2.687  | 3.2477 |
| 2.4568 | 4.8097 | 2.5512 | 2.3578 | 2.2786 | 3.8226 |
| 2.3957 | 2.3471 | 3.3985 | 2.7458 | 2.101  | 3.017  |

Table 4.3.2 Observed values of the objective function for 30 executions of simulated annealing, in the second economy

| 1.9004 | 1.6967 | 2.0286 | 2.7959 | 2.0924 | 0.0021 |
|---|---|---|---|---|---|
| 2.0924 | 2.7907 | 2.8026 | 2.0907 | 7.051  | 1.6967 |
| 1.8569 | 2.0286 | 7.051  | 5.4708 | 1.8682 | 3.4722 |
| 4.3208 | 5.2833 | 1.7628 | 2.3775 | 2.0384 | 2.2104 |
| 1.6967 | 1.8569 | 2.0286 | 2.7959 | 2.1983 | 2.0924 |

Table 4.3.3 Observed values of the objective function for 30 executions of stochastic local search, in the first economy

| 5.0852 | 3.4588 | 3.1157 | 3.0639 | 5.0852 | 3.6906 |
|--------|--------|--------|--------|--------|--------|
| 4.316  | 5.2471 | 4.4501 | 5.1925 | 4.43   | 4.8509 |
| 4.13   | 4.8509 | 5.7621 | 5.2317 | 5.3736 | 3.7022 |
| 3.7022 | 4.4699 | 3.5255 | 3.0277 | 5.1606 | 3.8351 |
| 4.5032 | 3.5166 | 5.366  | 4.7174 | 5.2953 | 5.7765 |

Table 4.3.4 Observed values of the objective function for 30 executions of stochastic local search, in the second economy

| 3.0554 | 4.3269 | 3.4749 | 4.6455 | 4.0029 | 3.8616 |
|--------|--------|--------|--------|--------|--------|
| 4.223  | 4.4017 | 3.1258 | 5.7564 | 4.7702 | 4.0098 |
| 4.6175 | 3.8385 | 4.7994 | 5.5739 | 4.3151 | 3.5352 |
| 4.7357 | 3.607  | 4.5146 | 3.3314 | 4.5831 | 4.0565 |
| 3.4749 | 5.5739 | 4.9971 | 3.8616 | 4.8727 | 5.0789 |

## 5. DISCUSSION

We have seen that the evolutionary algorithm of Ichibuchi et al., behaves in our expanded version of the model, in the same near-optimal way, as it did in the original model of Ichibuchi et al (2001). The differences in the actual numbers of the payoff function are to be expected, since the different placement of the firms and the markets, leads to different distances between each other, and consequently, to different actual profits. Therefore, the similarities displayed in figure 4.1 are enough to support the qualitative argument of similar behaviour between the two cases.

The significant differences in the observed values of the objective function (16), prove, in a statistical clear way, that different taxation coefficient yield different equilibrium situations in the economies. The mean value in the case of no taxes in the first economy is, for example 4.0108, while in the second case of table 4.2 it is 6.3299. The difference is so great, that we don't even need to perform a statistic hypothesis test, to draw the conclusion that these two samples come from different populations.

The same holds true for the differences in the relative performance of the two algorithms. Simulated annealing clearly outperforms the stochastic version of the local search algorithm we have introduced.

Since our executions are independent and we have used a sample consisting of 30 executions, we can assume that the objective function values are derived from a normal population (Koliva and Bora 1998). This allows us to evaluate probability intervals for the expected values of the algorithms. For stochastic search we have

| E (RD) | Market 1 | Market 2 |
|---|---|---|
| Sample mean | 4.464 | 4.301 |
| Sample deviation | 0.831 | 0.718 |
| max p.i. 95% | 4.762 | 4.557 |
| min p.i. 95% | 4.167 | 4.044 |
| max p.i. 98% | 4.817 | 4.605 |
| min p.i 98% | 4.112 | 3.996 |

So there is a 0.95 probability, for example, that the actual expected value for the objective function, that we would observe after 200 executions of the stochastic local search algorithm, in the first market, belongs in the interval
$$(4.167, 4.762)$$
For simulated annealing we have

| E (Bann) | Market 1 | Market 2 |
|---|---|---|
| Sample mean | 2.912 | 2.715 |
| Sample deviation | 0.916 | 1.589 |
| max p.i 95% | 3.239 | 3.284 |
| min p.i 95% | 2.584 | 2.146 |
| max p.i 98% | 3.301 | 3.390 |
| min p.i 98% | 2.522 | 2.040 |

The differences are great enough to prove in a statistical sound sense, that simulated annealing outperforms stochastic local search. The hypothesis test

$H_o: \mu_{SS} > \mu_{SA}$
$H_1: \mu_{SS} \leq \mu_{SA}$

is in favour of $H_1$, even with significance level (probability to reject $H_0$ while it is true) $a = 0.001$.

It is also apparent, from the above results, that the solution "proposed" by the simulated annealing algorithm, is far better, in terms of the objective function's value, than the case of no taxes at all.

## 6. CONCLUSIONS

We have successfully displayed that a meta-heuristic algorithm, simulated annealing, can be used in conjunction with simulation based on evolutionary algorithms, to provide a near-optimum decision of a peripheral economic planner, and help him in his quest to achieve a given goal. Although the model we have proposed is based on artificial data, since both the maps of the economies are totally artificial and the cost and demand functions provided by the model of Ichibuchi et al, are rather simplistic, it could be possible, for someone to try to access the same problem in a real economic context. This would require, of course, the correct mapping of real firms and markets in the model, plus some econometrically sound functions for the demand and cost functions. The evolutionary algorithm concerning the decisions of the agent should also be taken in consideration. This idea could hopefully be fruitful, both for future research and applications.